\documentclass[10pt]{article}
\usepackage{latexsym}
\usepackage{amssymb}
\usepackage{amsthm,amsmath}
\usepackage
          {graphicx}
\newtheorem{lemma}{Lemma}
\newtheorem{theorem}[lemma]{Theorem}

\newtheorem{corollary}[lemma]{Corollary}

\newtheorem{definition}[lemma]{Definition}

\newcommand{\E}{\mathbb E}
\newcommand{\PP}{\mathbb P}

\newcommand{\R}{\mathbb R}
\newcommand{\eps}{\varepsilon}

\title{Moment estimates for convex measures}

\author{
Rados{\l}aw Adamczak
 \thanks{Research partially supported by MNiSW Grant no. N N201 397437.}
\and Olivier Gu\'edon
\thanks{Research partially supported by ANR GeMeCoD, ANR 2011 BS01 007 01.}
\and  Rafa{\l} Lata{\l}a${}^{*}$
 \and  Alexander E. Litvak
\thanks{Research partially supported by  the
E.W.R. Steacie Memorial Fellowship.}
\and Krzysztof Oleszkiewicz${}^{*}$
\and Alain Pajor
\thanks{Research partially supported by the ANR project
ANR-08-BLAN-0311-01.}
\and  Nicole  Tomczak-Jaegermann
\thanks{This author holds the Canada Research Chair in
  Geometric Analysis.}
}

\newcommand\address{\noindent\leavevmode%
\noindent
Rados{\l}aw  Adamczak, \\
Institute of Mathematics, \\
University of Warsaw, \\
Banacha 2, 02-097 Warszawa, Poland\\
 \texttt{\small
e-mail:  radamcz@mimuw.edu.pl}

\medskip
\noindent
Olivier Gu\'edon,  \\
Universit\'{e} Paris-Est\\
\'{E}quipe d'Analyse et Math\'{e}matiques Appliqu\'ees, \\
5, boulevard Descartes,
Champs sur Marne,\\
77454 Marne-la-Vall\'{e}e,  Cedex 2, France\\
\texttt{\small%
e-mail: olivier.guedon@univ-mlv.fr}

\medskip
\noindent
Rafa{\l} Lata{\l}a, \\
Institute of Mathematics, \\
University of Warsaw, \\
Banacha 2, 02-097 Warszawa, Poland\\
 \texttt{\small
e-mail:   rlatala@mimuw.edu.pl}

\medskip
\noindent
Alexander E. Litvak, \\
Dept.~of Math.~and Stat.~Sciences,\\
University of Alberta, \\
Edmonton, Alberta, Canada, T6G 2G1.\\
\texttt{\small%
e-mail:   alexandr@math.ualberta.ca}

\medskip
\noindent
Krzysztof Oleszkiewicz, \\
Institute of Mathematics, \\
University of Warsaw, \\
Banacha 2, 02-097 Warszawa, Poland\\
 \texttt{\small
e-mail: koles@mimuw.edu.pl}

\medskip
\noindent
Alain  Pajor, \\
Universit\'{e} Paris-Est\\
\'{E}quipe d'Analyse et Math\'{e}matiques Appliqu\'ees, \\
5, boulevard Descartes,
Champs sur Marne,\\
77454 Marne-la-Vall\'{e}e,  Cedex 2, France\\
\texttt{\small%
e-mail: Alain.Pajor@univ-mlv.fr }

\medskip
\noindent
Nicole  Tomczak-Jaegermann, \\
Dept.~of Math.~and Stat.~Sciences,\\
University of Alberta, \\
Edmonton, Alberta, Canada, T6G 2G1.\\
\texttt{\small%
e-mail:    nicole.tomczak@ualberta.ca}

}

\begin{document}
\maketitle

\begin{abstract}
Let $p\geq 1$, $\eps >0$,  $r\geq (1+\eps) p$, and $X$ be a $(-1/r)$-concave
random vector in $\R^n$ with Euclidean norm $|X|$.
We prove that
$$
 (\E |X|^{p})^{1/{p}}\leq  c \left( C(\eps) \E|X|+\sigma_{p}(X)\right) ,
$$
where
$$\sigma_{p}(X)=\sup_{|z|\leq 1}(\E|\langle z,X\rangle|^{p})^{1/p}, $$
$C(\eps)$ depends only on $\eps$ and $c$ is a universal constant.
Moreover, if in addition $X$ is  centered then
$$
 (\E |X|^{-p})^{-1/{p}}\geq  c(\eps) \left(  \E|X| - C \sigma_{p}(X)\right) .
$$
\end{abstract}

\bigskip

\noindent {\bf 2010 Math. Subject Classification}:  Primary: 46B06, 60E15, 60F10;
\newline  Secondary: 52A23, 52A40  \\

\noindent {\bf Key Words and Phrases:} convex measures,
$\kappa$-concave measure, tail inequalities, small ball probability estimate.

\section{Introduction}

Let $X$ be a random vector with values in a finite dimensional Euclidean space $E$
with Euclidean norm $|\cdot|$ and scalar product  $\langle \cdot,\cdot\rangle$.
 For any $p>0$, we define
the weak $p$-th moment of $X$  by
\[
\sigma_{p}(X)=\sup_{|z|\leq 1}(\E
|\langle z,X\rangle|^p)^{1/p}.
\]

Clearly $ (\E |X|^p)^{1/p} \geq
\sigma_{p}(X)$ and by  H\"{o}lder's inequality,
 $ (\E |X|^p)^{1/p} \geq \E |X| $. In this paper we are interested in reversed inequalities
of  the form
\begin{equation}
\label{equa:relationship}
 (\E |X|^p)^{1/p} \leq  C_1\, \E|X|   +   C_2 \sigma_p(X)
\end{equation}
for  $p\geq 1$ and constants $C_1$ and $C_2$.

This is known for some classes of distributions and the question has been studied in a more
general setting (see \cite{La} and references therein) and our objective in this paper
is to describe new classes for which the relationship (\ref{equa:relationship}) is satisfied.

Let us recall some known results when (\ref{equa:relationship}) holds.
It clearly holds for Gaussian vectors and it is not
difficult to see that (\ref{equa:relationship}) is true for subgaussian vectors
(see below for definitions) for every $p\geq 1$, with $C_1$ and $C_2$ depending
only on the subgaussian parameter.

Another example of such a class is the class of  so-called log-concave vectors.
A probability measure $\mu$ on $\R^m$ is called log-concave
if for all $0<\theta<1$ and
for all compact subsets $A,B\subset\R^m$ with positive measure one has
\begin{equation}\label{deflog}
\mu( (1-\theta) A+\theta B)\geq
\mu(A)^{1-\theta}\mu(B)^\theta.
\end{equation}
A random vector with a log-concave distribution is called
log-concave.
It is  known that
for every log-concave random vector $X$ in a finite dimensional Euclidean space and any $p>0$,
\[
(\E|X|^p)^{1/p}\leq C\big(\E |X|+\sigma_{p}(X)\big),
\]
where $C>0$ is a universal constant.
See Corollary \ref{corollary:imprPaouris} and references below.

In this paper we  will consider the class of convex measures introduced
by Borell.  Let $\kappa<0$.
A probability measure $\mu$ on $\R^m$ is called
$\kappa$-concave if for all $0<\theta<1$ and for all compact
subsets $A,B\subset\R^m$ with positive measure one has
\begin{equation}\label{defconc}
\mu( (1-\theta) A+\theta B)\geq
\left((1-\theta)\mu(A)^\kappa +\theta \mu(B)^\kappa\right)^{1/\kappa}.
\end{equation}
A random vector with a $\kappa$-concave distribution is called
$\kappa$-concave. Note that a log-concave vector is also $\kappa$-concave
for any $\kappa<0$.

\smallskip

We show in  Theorem \ref{theorem:convexPaouris}  that  for $\kappa > -1$, a
$\kappa$-concave random vector satisfies (\ref{equa:relationship})
for all $0<(1+\eps) p<-1/\kappa$ with $C_1$ and $C_2$ depending only on $\eps$.

\smallskip

In fact, in Definition~\ref{lemma:ratio logconcave} we will introduce a general
assumption on the distribution, called $H(p,\lambda)$.  The main result of the
first part of the paper is  Theorem \ref{theorem:imprPaouris2}  in which we
show that this assumption is  sufficient in order to have (\ref{equa:relationship}).
In Theorem~\ref{theorem: assumption  convex} we prove that convex measures satisfy
this assumption.

\smallskip

One of the main applications  of  the relationship  (\ref{equa:relationship})
consists in tail inequalities  for $ \PP \left( |X|\ge t\, \E|X|\right)$.
In  Corollary~\ref{corollary:convex-deviation-1} we show that for $r>2$ and
for a $(-1/ r)$-concave isotropic random vector $X\in\R^n$ the above probability
is bounded by
$\left(\frac{c\max\{1, r/\sqrt n\}}{t}\right)^{ r/2}$.
From this bound we deduce that
the empirical  covariance matrix of a sample of size proportional to $n$ is a good
approximation of the covariance matrix of  $X$, extending results of
\cite{ALPTjams, ALPTcras} from log-concave measures to convex measures.
This provides thus a new class of random vectors satisfying such approximation.
See Corollary~\ref{corollary:SV} and the remark following it.

\smallskip

The second part of the paper deals with negative moments. We are looking for relationship of
the  form
\begin{equation}
\label{equa:relationship2}
   \left( \E |X|^{-p}\right)^{-1/p} \geq  C_1 \,  \E |X| - C_2 \sigma_p(X)
\end{equation}
for  $p>0$ and constants $C_1$ and $C_2$.

We show in Theorem~\ref{theorem:negative moment} that  for $\kappa > -1$, an
$n$-dimensional  $\kappa$-concave random vector satisfies (\ref{equa:relationship2})
for all $0<(1+\eps) p<\min\{n/2,(-1/\kappa)\}$ with $C_1$ and $C_2$ depending only on $\eps$.
As an application we show a small ball probability estimate for $\kappa$-concave random vectors.
In  the log-concave setting it was proved in \cite{Pa2}.

\section{Preliminaries}

The space $\R^m$ is equipped with the scalar product $\langle
\cdot,\cdot\rangle$, the Euclidean norm $|\cdot|$, the unit ball
$B_2^m$ and the volume measure $\text{vol}(\cdot)$. The canonical basis
is denoted by $e_1$, $e_2$, \ldots, $e_m$. A gauge or
Minkowski functional $\|\cdot \|$ on $\R^m$ is a non-negative function
on $\R^m$ satisfying: $\| \lambda x \|=\lambda \| x\|$ and $\|
x+y\|\leq \| x\|+\|y\|$ for every $x,y\in\R^m$ and every real
$\lambda\geq 0$ and such that $\| x\|$=0 if and only if $x=0$. The
dual gauge is defined for every $x\in\R^m$ by $\| x\|^*=\max\{
\langle x,t\rangle: \|t\|\leq 1\}$.  A  body is a compact
subset of $\R^m$ with a non-empty interior.  Any convex body $K\subset \R^m$
containing the origin in its interior defines the gauge by
$\|x\|=\inf\{ \lambda\geq 0\,:\, x\in \lambda K\}$. It is called the
Minkowski functional of $K$.  If $K\subset \R^m$ is a convex body
containing the origin in its interior, the polar body $K^\circ$ is
defined by  $K^\circ=\{x\in\R^m: \langle x,t\rangle\leq 1 \ \text{for
 all}\ t\in K\}$.  The diameter of $K$ in the Euclidean metric is
denoted by $\text{diam} (K)$.

{}For a linear subspace $F\subset \R^n$ we denote the orthogonal
projection on $F$ by $P_F$. Note that $P_F K^{\circ} := P_F (K^{\circ})
 = (K\cap F)^{\circ}$, when the polar is taken in $F$.

{}For a random vector $X$ in $\R^n$ with a density $g$ and a subspace
$F\subset \R^n$, we denote the density of $P_F X$ by $g_F$.

A random vector $X$ in $\R^m$ will be called non-degenerated if
it is not supported in a proper affine subspace of $\R^m$.
It is called isotropic if it is centered and
for all $\theta\in\R^m$, $\E|\langle X,\theta\rangle|^2=|\theta|^2$.

Given a non-negative bounded function $g$ on $\R^m$ we introduce the
following associated set. For any $\alpha\geq 1$, let
\begin{equation}\label{kalpha}
K_{\alpha}(g)=\{t\in\R^m: g(t)\geq \alpha^{-m}\|g\|_\infty\},
\end{equation}
where $\|g\|_\infty=\sup_{t\in\R^m} |g(t)|$.

By $g_i$, $g_{i, j}$ we denote independent standard Gaussian
random variables, i.e. centered and of variance one. A standard
Gaussian vector in $\R^n$ is denoted by $G$, i.e. $G=(g_1, g_2, ...,
g_n)$. The standard Gaussian matrix is the matrix whose entries are
i.i.d. standard Gaussian variables, i.e. $\Gamma = \{g_{i,j}\}$. By
$\gamma_p$ we denote the $L_p$ norm of $g_1$.  Note that $\gamma_p
/\sqrt{p} \to 1/\sqrt{e}$ as $p\to\infty$.

We denote by $\mu _{n,k}$ the  Haar probability measure on the Grassmannian
$G_{n,k}$ of $k$-dimensional subspaces of $\R^n$.

Recall that for a real number $s$, $\lceil s\rceil$ denotes the smallest integer
which is not less than $s$.

By $C$, $C_0$, $C_1$, $C_2$, ..., $c$, $c_0$ $c_1$, $c_2$ we denote
absolute positive constants, whose values can change from line to
line.

For two functions $f$ and $g$ we write $f\sim g$ if there are absolute positive
constants $c$ and $C$ such that $c f\leq g\leq C f$.

\section{Convex probabilities}

Let $\kappa\leq 1/m$. A Borel probability measure $\mu$ on $\R^m$ is called
$\kappa$-concave if it satisfies (\ref{defconc}). When $\kappa=0$, this inequality
should be read as (\ref{deflog}) and it defines $\mu$ as a log-concave probability.

In this paper we will be interested in the case $\kappa\leq 0$, which we consider
from now on.

The class of $\kappa$-concave measures was introduced and studied by
Borell.  We refer to \cite{Bo1,Bo2} for a general study and to
\cite{BobM} for more recent development.  A $\kappa$-concave
probability is supported on some convex subset of an affine subspace
where it has a density.  When the support generates the
whole space,  a characterization of Borell (\cite{Bo1,Bo2}) states
that the probability is absolutely continuous with respect to the Lebesgue
measure and has a density $g$ which is log-concave when $\kappa=0$ and
when $\kappa<0$, is of the form
\[ g=f^{-\beta}\quad\text{with}\quad \beta =m-\frac{1}{\kappa}, \]
where $f:\R^m \to (0, \infty ]$ is a convex function.  The class of
$m$-dimensional $\kappa$-concave probabilities is increasing as $\kappa$
is decreasing. In particular a log-concave probability is
$\kappa$-concave for any $\kappa<0$.

As we mentioned in the Introduction, a random vector with a $\kappa$-concave
distribution is called $\kappa$-concave.  Clearly, the linear image of a $\kappa$-concave
probability is also $\kappa$-concave.  Recall that any semi-norm
of an $m$-dimensional vector with a $\kappa$-concave distribution has
moments up to the order $p<-1/\kappa$ (see \cite{Bo1} and Lemmas~\ref{lemma:khinchin}
and \ref{medianone} below). Since we are interested in comparison of moments
with the moment of order 1, we will always assume that $-1<\kappa\leq 0$.

\section{Strong and weak moments}

In this section we consider a random vector $X$ in a finite
dimensional Euclidean space $E$.

\begin{definition}
\label{lemma:ratio logconcave}
Let $p>0$, $m=\lceil p\rceil$, and $\lambda\geq 1$. We say that a
random vector $X$ in $E$ satisfies the assumption $H(p,\lambda)$
if for
every linear mapping $A:E\to\R^m$ such that $Y=AX$ is
non-degenerate
there exists a gauge $\|\cdot\|$ on $\R^m$ such that
$\E\|Y\|< \infty$ and
\begin{equation}\label{equ:hypthesis H1}
  (\E\|Y\|^p)^{1/p}\leq \lambda\,\E \|Y\|.
\end{equation}
\end{definition}

\medskip

\noindent {\bf Remark.}
Let us give a  first example of a random
vector satisfying $H(p,\lambda)$. Let $X$ be a random vector in an
$n$-dimensional Euclidean space $E$, satisfying,
for some  $\psi \geq 1$,
\begin{equation}
\label{equ:subgaussian}
  \forall z\in E \quad\forall\, 1\leq p\leq n\quad (\E |\langle
  z,X\rangle|^p)^{1/p}\leq \psi \sqrt p\,
  \E |\langle z,X\rangle|.
\end{equation}
Then $X$ satisfies $H(p, C \psi^2)$ for every  $p\geq 1$.
For example, the standard Gaussian and
Rademacher vectors  satisfy the above condition with
$\psi$ being a numerical
constant.  More generally, if $X$ is subgaussian,
then $X$ satisfies (\ref{equ:subgaussian}).

To prove that $(\ref{equ:subgaussian})$ implies $H(p, C \psi^2)$,
let $p>0$, $m=\lceil p\rceil$ and let  $A:E\to \R^m$  be  such that
$Y=AX$ is non-degenerate. We may assume that $m\geq 2$. Clearly,
because of the linear invariance of the
property~(\ref{equ:subgaussian}), we may also assume that $Y=AX$ is
isotropic.
Thus
$(\ref{equ:subgaussian})$
yields,
\begin{equation}
\label{equ:hypothesis H3}
\forall  z\in \R^m \quad
(\E |\langle z,Y\rangle|^p)^{1/p}
\leq \psi
  \sqrt p\,
  \E |\langle z,Y\rangle|\leq \psi\sqrt m |z|
\leq \sqrt{2} \psi^2 |z|    \E |Y|,
\end{equation}
where the last inequality follows from isotropicity of $Y$ by applying
$(\ref{equ:subgaussian})$ with $p=2$, $z_i=A^\top e_i$, $i\leq m$,
and the Cauchy-Schwarz inequality.

Now let us make the following general observation. Let $p\geq 1$ and
$m=\lceil p\rceil$.  Let $Y$ be a
random vector in an $m$-dimensional normed space with norm
$\|\cdot\|$.  Since any $m$-dimensional norm can be estimated,
up to a multiplicative constant, by the supremum over an
exponential  (in $m$)
number of norm one linear forms, we deduce that
\begin{equation}
\label{equ:hypothesis H2}
   \left(\E \|Y\|^p\right)^{1/p}\leq C'
           \sup_{\|\varphi\|^* \leq 1}
   \left(\E|\varphi(Y)|^p\right)^{1/p},
\end{equation}
where $C'$ is a universal constant (see \cite{LW}
Proposition 3.20).
Combining this  with
(\ref{equ:hypothesis H3})

we conclude that
$$
\left(\E |Y|^p\right)^{1/p}\leq C' \sup_{|z| \leq 1}
      (\E |\langle z,Y\rangle|^p)^{1/p}
\leq C\,C'\,\psi\, \E |Y|,
$$
which shows that $X$ satisfies $H(p, C C'\psi)$.

\medskip

The main result of this section states a relationship between weak
and strong moments under the assumption $H(p,\lambda)$.

\begin{theorem}
\label{theorem:imprPaouris2}
Let $p>0$ and $\lambda\geq 1$. If  a random vector $X$  in a finite dimensional Euclidean space
satisfies $H(p,\lambda)$, then
$$
  (\E |X|^p)^{1/p} \leq  c(\lambda\, \E|X|   +    \sigma_p(X) ),
$$
where $c$ is a universal constant.
\end{theorem}

The first step of the proof of Theorem \ref{theorem:imprPaouris2} consists of showing that
there exists some $z$  such that $(\E(\langle z, Y\rangle)_+^{p})^{1/p }$ is small, with
comparison to $\E |Y|$.  This is the purpose of the following lemma.

\begin{lemma}
\label{lemma:minimum}
Let $Y$ be a random vector in $\R^m$. Let $\|\,\cdot\,\|_1$ and
$\|\,\cdot\,\|_2$ be two gauges on $\R^m$ and  $\|\,\cdot\,\|_1^*$
and $\|\,\cdot\,\|_2^*$ be their dual gauges.
 Then for all $p >0$,
\[
\min_{\|z\|_2^*=1}
(\E(\langle z, Y\rangle)_+^{p})^{1/p } \leq
\frac{(\E\|Y\|_1^{p})^{1/p}}{\E\|Y\|_1 }\,
\E\|Y\|_2.
\]
\end{lemma}

\noindent{\bf Proof. }
Let  $r$ be the largest real number such that
$r\|t\|_1 \le \|t\|_2$  for all $t \in \R^m$.
By duality  $r$ is  the largest number such that
$r\|w\|_2^* \le \|w\|_1^*$  for  $w \in \R^m$.
Pick  $z\in \R^m$ such that  $\|z\|_2^*=1$
and $\|z\|_1^* = r$. Then
$\langle z, t \rangle \leq \|z\|_1^*\|t\|_1 \leq r \|t\|_1$
for all  $t\in \R^m$.
Therefore, for any $p>0$,
$(\E(\langle z, Y \rangle)_+^{p})^{1/p}
\leq r\,(\E\|Y\|_1^{p})^{1/p}$. Thus  the lemma  follows from the
inequality
$r \,\E\|Y\|_1\leq \E\|Y\|_2 $. \qed

\medskip
The second  step of the proof of Theorem \ref{theorem:imprPaouris2}
is contained in the next lemma.

\begin{lemma}
\label{lemma:minimax}
Let $n, m\geq 1$ be integers. Let $p\geq 1$. Let $X$ be an
$n$-dimensional random vector and $\Gamma$ be an $n\times m$
standard Gaussian matrix. Then
\[
   (\E |X|^p)^{1/p}\leq 2^{1/p}\, \gamma_p^{-1} \left(\E\min_{|t|=1}
   \|\Gamma t\|_{p,+} + (C\gamma_p+ \sqrt m)\, \sigma_p(X)\right),
\]
where $\|z\|_{p,+}=\left(\E\left(\langle z,
X\rangle\right)_+^{p}\right)^{1/p}$
and $C$ is a universal constant.
\end{lemma}

\noindent {\bf Proof. } For every $x,y\in\R^n$, $\left|\, \|x\|_{p,+}
- \,  \|y\|_{p,+}\,\right| \leq |x-y|\sigma_p(X).$
The classical Gaussian
concentration inequality (see \cite{CIS} or inequality (2.35) and
Proposition~2.18 in \cite{L})  gives that
\[
\PP \left(\bigl|\, \|G\|_{p,+} - \E \|G\|_{p,+}\,\bigr|\geq s\right) \leq
2 \exp \left({-s^2/2 \sigma_p^2(X)}\right) ,
\]
and implies (cf. \cite{LMS}, Statement 3.1)
\begin{equation}
\label{equ:gaussian moment}
(\E\| G\|_{p,+}^p)^{1/p} \leq \E\| G\|_{p,+}+C\gamma_p\sigma_p(X),
\end{equation}
where $C$ is a universal constant.
Since $\langle G, X\rangle$ has the same distribution as $|X|\, g_1$,

we have
\begin{equation}
 \E \left(\langle G, X\rangle\right)_+^{p} =
   ({1}/{2})\,
 \E|\langle G, X\rangle|^{p}
 \quad \mbox{ and } \quad
 \E|\langle G, X\rangle|^{p}=\gamma_p^p \,\E |X|^p .
 \end{equation}
Therefore
$$
  \left(\E |X|^p \right)^{1/p} = 2^{1/p}\, \gamma_p^{-1} \left(\E\|G\|_{p,+}^p
  \right)^{1/p}  \leq 2^{1/p}\, \gamma_p^{-1}\big(\E\|G\|_{p,+}+
   C\gamma_p\sigma_p(X)\big).
$$
The  Gordon minimax lower bound (see \cite{Go},  Theorem 2.5)
states that for any norm $\| \cdot  \|$
\[
  \E \|G\|\leq \E\min_{|t|=1}\|\Gamma t\|+ \E |H|\,\max_{|z|=1}\|z\|,
\]
where  $H$ is a  standard Gaussian vector in $\R^m$.
It is easy to check  the proof and to show that this inequality
remains true  when $\| \cdot  \|$
is a gauge. This gives us that
\[
  \E \| G\|_{p,+}\leq \E\min_{|t|=1}\|\Gamma t\|_{p,+}+ \E |H|\,\max_{|z|=1}\|z\|_{p,+}
  \leq \E\min_{|t|=1}\|\Gamma t\|_{p,+}+\sqrt m\max_{|z|=1}\|z\|_{p,+}
\]
and it is enough to observe that $\max_{|z|=1}\|z\|_{p,+}
\leq\sigma_p(X)$.
\qed

\medskip

\noindent{\bf Proof  of  Theorem \ref{theorem:imprPaouris2}.}
We may assume that $p\geq 1$ since  Theorem~\ref{theorem:imprPaouris2}
is obviously true when $0<p\leq 1$.
Let $m$ be the integer so that $1\leq p\leq m<p+1$, thus $m\leq 2p$.
 We  use the notation of  Lemma \ref{lemma:minimax}.
  We first condition on $\Gamma$.
We have
\[\| \Gamma z\|_{p,+}=(\E_X(\langle \Gamma z, X\rangle)_+^{p})^{1/p}=
(\E_X(\langle z, \Gamma^*X\rangle)_+^{p})^{1/p}.\]
Let $Y=\Gamma^*X\in \R^m$. If $Y$ is supported by a hyperplane then
$$\min_{|z|=1}(\E_X(\langle z, \Gamma^*X\rangle)_+^{p})^{1/p}=0.$$
Otherwise by our assumption $H(p,\lambda)$
there exists a gauge in $\R^m$ such that
$$
    (\E\|Y\|^p)^{1/p}\leq \lambda\,\E \|Y\|.
$$
{}From Lemma \ref{lemma:minimum} we get
\[
    \min_{|z|=1}(\E_X(\langle z, \Gamma^*X\rangle)_+^{p})^{1/p}\leq
    \lambda\E_X|\Gamma^*X|.
\]
We now take the expectation with respect to $\Gamma$ and get
\[
  \E\min_{|t|=1}\| \Gamma t\|_{p,+} \leq  \lambda\E|\Gamma^*X|
  =\lambda\E |H|\E|X|\leq \lambda\,\sqrt m \, \E|X|,
\]
where $H$ is a standard $m$-dimensional Gaussian vector.
The proof is concluded using  Lemma \ref{lemma:minimax}  and the fact
that
$\gamma_p^{-1}\sqrt p$ is bounded. Indeed,
\[
  \begin{array}{rl}
    (\E |X|^p)^{1/p}&\leq 2^{1/p}\gamma_p^{-1} \lambda\,\sqrt m\,
    \E|X| + 2^{1/p}
    (C+\gamma_p^{-1}\sqrt m) \sigma_p(X)\\
    &\leq c'\left(\lambda\E|X|+  \sigma_p(X)\right).
 \end{array}
\]\qed

\section{Tail behavior of convex measures}

\begin{theorem}
\label{theorem: assumption  convex}
Let $n\geq 1$ and $r>1$. Let $X$ be a centered $(-1/r)$-concave
random vector in a finite dimensional Euclidean space. Then
for every $0<p<r$, $X$ satisfies the assumption $H(p, \lambda(p,r))$
with $\lambda(p,r)= c\left(\frac{r}{r-1}\right)^3\left(\frac{r}{r-p}\right)^4$,
where $c$ is a universal constant.
\end{theorem}

\noindent {\bf Remark:}  Note that the parameter $\lambda(p,r)$ in
Theorem~\ref{theorem: assumption  convex} is bounded
by a universal constant if the parameters $p$ and $r$ are not close,
for instance if $r\geq 2\max\{1, p\}.$

\bigskip

\begin{theorem}
\label{theorem:convexPaouris}
Let $r> 1$ and let $X$ be a $(-1/ r)$-concave
random vector in a finite dimensional Euclidean space. Then, for every    $0<p<r$,
\begin{equation}\label{thsix}
  (\E |X|^{p})^{1/{p}}\leq c(C_2(p,r) \E|X|+ \sigma_p(X)),
\end{equation}
where $C_2(p,r)=c\left(\frac{r}{r-1}\right)^3\left(\frac{r}{r-p}\right)^4$
and $c$ is a universal constant.
\end{theorem}

\noindent {\bf Proof.} The proof may be reduced to the case of a centered random vector.
Indeed, let $X$ be a $(-1/ r)$-concave random vector, then so is $X - \E X$.
Since
$$
(\E|X|^p)^{1/p} \le (\E|X - \E X|^p)^{1/p} + |\E X|\le (\E|X - \E X|^p)^{1/p} + \E |X|,
$$
$\E |X-\E X|\leq 2 \E |X|$ and
$\sigma_p(X-\E X)\leq 2\sigma_p(X)$,
 we may assume that $X$ is centered.
The  theorem now follows immediately by Theorems \ref{theorem:imprPaouris2}
and \ref{theorem: assumption  convex}.
\qed

\medskip

Note that  trivially a reverse inequality to (\ref{thsix}) is valid, for $p\ge 1$:
\[
2 (\E |X|^{p})^{1/{p}}\geq
\E|X|+ \sigma_{p}(X).
\]
Therefore Theorem \ref{theorem:convexPaouris} states an equivalence
\[
(\E |X|^{p})^{1/{p}} \sim_{C_2(p,r)}
\E|X|+ \sigma_{p}(X).
\]

Since a log-concave measure is $\kappa$-concave for any $\kappa< 0$,
we obtain

\begin{corollary}
\label{corollary:imprPaouris}
For any log-concave random vector $X$ in a finite dimensional Euclidean space and any $p>0$,
\[
(\E|X|^p)^{1/p}\leq C\big(\E |X|+\sigma_{p}(X)\big),
\]
where $C>0$ is a universal constant.
\end{corollary}

As formulated here, Corollary \ref{corollary:imprPaouris}
first appeared as Theorem 2 in \cite{cras_allpt} (see also \cite{ALLPT}).
A short proof of this result was given in \cite{ALLOPT}. It can be deduced
directly from Paouris work in \cite{Pa} (see \cite{ALLOPT}).

\medskip

As it was mentioned above, if $X\in E$ is $(-1/r)$-concave then so is
$\langle z,X\rangle$ for any $z\in E$. From Lemma
\ref{lemma:khinchin}, we have that for any $1\leq p<r$, \begin{equation}\label{equ:weak moment inequ}
(\E|\langle z,X\rangle|^p)^{1/p}\leq C_1(p,r)\, \E|\langle
z,X\rangle| ,
\end{equation}
where $C_1(p,r)$ is defined in Lemma~\ref{lemma:khinchin}.

Assume that $r>2$. Let $n$ be the dimension of $E$. If moreover $X$ is centered and has the identity as the
covariance matrix -- such a vector is called an {\em isotropic random
  vector} -- then one has for any $z\in S^{n-1}$ and any $1\leq p<r$,
\begin{equation}\label{logiso}
  (\E|\langle z,X\rangle|^p)^{1/p}\leq C_1(p,r)\, \E|\langle
  z,X\rangle|\leq C_1(p,r)  (\E|\langle z,X\rangle|^2)^{1/2}=C_1(p,r).
\end{equation}

Since in that case, $\E |X|\le (\E |X|^2)^{1/2}=\sqrt n$, it follows  from Theorem
\ref{theorem:convexPaouris}
that for any $1\leq p<r$,
\begin{equation}\label{forcor}
\left(\E |X|^{p}\right)^{1/{p}}\leq c(C_2(p,r)\sqrt n+ C_1(p,r)).
\end{equation}

Together with Markov's inequality this give the following
Corollary.

\begin{corollary}
\label{corollary:convex-deviation-1}
Let $r> 2$ and  let  $X\in\R^n$ be a $(-1/ r)$-concave isotropic
random vector.  Then for every $t> 0$,
\begin{equation}\label{equ:paouris1}
  \PP \Big( |X|\ge t\sqrt n\Big)\leq
  \left(\frac{c\max\{1, r/\sqrt n\}}{t}\right)^{ r/2} .
\end{equation}
In particular, if $r\geq 2\sqrt n$, then for every $6c\leq t\leq 3cr/\sqrt n$,
\begin{equation}\label{equ:paouris2}
  \PP \Big( |X|\ge t\sqrt n\Big)\leq \exp(-c_0 t\sqrt n),
\end{equation}
where $c$ and $c_0$ are universal positive constants.
\end{corollary}

\medskip

\noindent
{\bf Remark. }
A log-concave measure is $(-1/r)$-concave for every $r>0$, thus in such
a case inequality (\ref{equ:paouris2}) is valid for every $t>c$,
which is a result from \cite{Pa}.

\medskip

\noindent {\bf Proof of Corollary~\ref{corollary:convex-deviation-1}. }
The inequality~(\ref{equ:paouris1}) follows by Markov's inequality
from inequality~(\ref{forcor}) with $p=r/2$, since $C_2(r/2, r)\leq c$ and
$C_1(r/2, r)\leq c r$ for a universal positive~$c$.

To prove the ``In particular" part denote $r'= t\sqrt n/(3c)$. Note that
 $r'\geq 2\sqrt n$ and that  $r' \leq r$. Therefore $X$ is $(-1/r')$-concave
as well and we can apply (\ref{equ:paouris1}) with $r'$, obtaining the bound
for probability $3^{-r'/2}$, which implies the result.
\qed

\bigskip

We now apply our results to the problem of the approximation of
the covariance matrix by the empirical covariance matrix. Recall that for
a random vector $X$ the covariance matrix of $X$ is given by $\E X X^\top$.
It is equal to the identity operator $I$ if $X$ is isotropic. The empirical
covariance matrix of a sample of size $N$ is defined by $\frac{1}{N}\sum_{i=1}^N  X_iX_i^\top$,
where $X_1, X_2, \dots, X_N$ are independent copies of $X$. The main
question is how small $N$ can be taken  in order that these two matrices are close
to each other in the operator norm (clearly, if $X$ is non-degenerated then
$N\geq n$ due to the dimensional restrictions and, by the law of large numbers,
the empirical covariance matrix tends to the covariance matrix as $N$ grows to infinity).
See \cite{ALPTjams, ALPTcras} for references on this question and for corresponding
results in the case of log-concave measures. In particular, it was proved there that
for $N\geq n$ and log-concave $n$-dimensional vectors $X_1,\cdots,X_N$ one has
$$
  \Big\| \frac{1}{N}\sum_{i=1}^N  X_iX_i^\top-I\Big\| \leq C \sqrt{\frac{n}{N}}
$$
with probability at least $1-2\exp(-c\sqrt{n})$,
where, as usual, $I$ is the identity operator, $\|\,\cdot\,\|$ is the operator norm
$\ell_2^n\to \ell_2^n$ and  $c$, $C$ are absolute positive constants.

In  \cite{SV} (Theorem~1.1), the following condition was introduced:
an isotropic random vector $X\in\R^n$ is said to satisfy the
{\em strong regularity assumption} if for some
$\eta,C>0$ and every rank $k\leq n$ orthogonal projection $P$, one has for every $t\geq C$
$$
  \PP\left(|PX| \geq t \sqrt{k} \right) \leq C \, t^{-2-2\eta} k^{-1-\eta}.
$$

We show that an isotropic $(-1/r)$-concave random vector satisfies this assumption.
For simplicity we will show this with $\eta =1$ (one can change $\eta$ adjusting
constants).

\begin{lemma}
 Let $n\geq 1$, $a>0$ and $r=\max\{4, 2a\log n\}$.
Let $X\in\R^n$ be an isotropic $(-1/r)$-concave random vector.
 Then there exists an absolute constant $C$ such that for
 every rank $k$ orthogonal projection $P$ and every
 $t\geq C_1(a)$, one has
$$
  \PP\left(|PX| \geq t \sqrt{k} \right) \leq C_2(a) \, t^{-4} k^{-2} ,
$$
where $C_1(a)=C  \exp{(4/a)}$ and $C_2(a)=C \max\{(a \log a)^4, \exp{(32/a)}\}$.
\end{lemma}

\medskip

\noindent
{\bf Proof. }
Let $P$ be a projection of rank $k$. Let $c$ be the constant from
Corollary~\ref{corollary:convex-deviation-1} (without loss of generality
we assume $c\geq1$) and $t>c$. If $r\leq  \sqrt{k}$
then Corollary~\ref{corollary:convex-deviation-1} implies that
$$
 \PP \Big( |PX|\ge t\sqrt k\Big)\leq  \left(\frac{c}{t}\right)^{r/2}
 \leq \left(\frac{c}{t}\right)^{a \, \log n} = n^{-a\log (t/c)}.
$$
If $r>\sqrt{k}$ then $X$ is also $(-1/\sqrt{k})$-concave, hence,  assuming
$k>\max\{4, 64 a^2 \log^2(4a)\}$ (so that $\sqrt{k}\geq 2a \log k$) and
applying Corollary~\ref{corollary:convex-deviation-1} again
we obtain
$$
 \PP \Big( |PX|\ge t\sqrt k\Big)\leq  \left(\frac{c}{t}\right)^{\sqrt{k}/2}
 \leq \left(\frac{c}{t}\right)^{a \, \log k} = k^{-a\log (t/c)}.
$$
Thus in both cases we have
$$
  \PP \Big( |PX|\ge t\sqrt k\Big)\leq   k^{-a\log (t/c)}.
$$
One can check that for $t\geq c^2 \exp{(4/a)}$ and $k\geq \exp{(16/a)}$
this implies
$$
  \PP\left(|PX| \geq t \sqrt{k} \right) \leq  t^{-4} k^{-2},
$$
which proves the desired result for $k> C_a:=\max\{64 a^2 \log^2(4a), \exp{(16/a)}\}$
and $t\geq c^2 \exp{(4/a)}$.

Assume now that $k\leq C_a$.
Then we apply Borell's Lemma -- Lemma~\ref{lemma:khinchin}.
We have that for every $t\geq 3$
$$
  \PP\left(|PX| \geq t \sqrt{k} \right) \leq
  \left(1+\frac{t}{9r}\right)^{-r}.
$$
It is not difficult to see (e.g., by considering cases $t\leq 9r$, $9r<t\leq 18r$ and
$t>18r$) that for $C(a) := 54^4 C_a^2$, $t\geq 3$ and $r\geq 4$, one has
$$
  \PP\left(|PX| \geq t \sqrt{k} \right) \leq C(a) t^{-4} k^{-2}.
$$
This completes the proof.
\qed

\medskip

Theorem~1.1 from  \cite{SV} and the above lemma immediately imply the
following corollary on the approximation of the covariance matrix by
the empirical covariance matrix.

\begin{corollary}
\label{corollary:SV}
 Let $n\geq 1$, $a>0$ and $r=\max\{4, 2a\log n\} $.  Let $X_1,\dots, X_N$ be
independent $(-1/ r)$-concave isotropic random vectors in $\R^n$.
Then for every $\varepsilon \in (0,1)$ and every $N\geq
C(\varepsilon,a) n$, one has
$$
   \E\Big\| \frac{1}{N}\sum_{i=1}^N  X_iX_i^\top-I\Big\| \leq\varepsilon,
$$
where
$C(\varepsilon,a)$ depends only on $a$ and $\varepsilon$.
\end{corollary}

\noindent{\bf Remark. } Let $r=2a\log(2n)>8$.  Applying
Corollary~\ref{corollary:convex-deviation-1} for independent $(-1/
r)$-concave isotropic random vectors $X_1$, $X_2$, \dots,$X_N$
and using  results of \cite{MendPao}, it can be checked that with large probability
$$
  \Big\| \frac{1}{N}\sum_{i=1}^N  X_iX_i^\top-I\Big\| \leq C(a)\sqrt{\frac{n}{N}}
$$
 where $C(a)$ depends only on $a$.
As we mentioned above, this extends the results of \cite{ALPTjams, ALPTcras}
on the approximation of the covariance matrix from the log-concave setting
to the class of convex measures.

\bigskip

Now we prove Theorem~\ref{theorem: assumption  convex}.
We need the following lemma. Recall that $K_{\alpha}$ was defined
by (\ref{kalpha}).

\begin{lemma}
\label{lemma: restriction}
Let $m$ be an integer. Let $r> 1$ and $0< p< r$. Let $Y\in\R^m$ be a
centered random vector with density $g=f^{-\beta}$
with  $\beta=m+r$ and $f$ convex positive.
Let $F:\R^m\to \R^+$ be such that for every $t\in\R^m$,
$F(2t)\leq 2^p F(t)$ and assume that $\E \/ F(Y)$ is finite.
Then, there exists a positive universal constant $c$ such that $0\in K_\alpha(g)$ and
\begin{equation}
\label{equ:restriction}
 \E F(Y)\leq c(p,r)\,\E \left(F(Y)1_{K_\alpha(g)} (Y)\right),
 \end{equation}
where $c(p,r)=1+\frac{c}{r-p}$ and
$\alpha=\left(c\frac{(m+r)^2}{(r-p)(r-1)}\right)^{\frac{m+r}{m}}$.
\end{lemma}

\noindent{\bf Proof. }
Let $\alpha>0$ be specified later.
Let $\gamma=\frac{\beta-m-1}{\beta-1}$. From Lemma  \ref{corollary:g(0)} we have $\gamma f(0)\leq
\min f=\|g\|_\infty ^{-1/\beta}$ and by definition,
 $\min f<\alpha^{-m/(r+m)}f(t)$ when $t\notin K_\alpha(g)$.
Using the convexity of $f$ and the last two inequalities we get
\begin{equation}\label{equ:convex1}
\forall t\notin K_\alpha(g)\quad  g(t/2)\geq g(t)\left( \frac{1}{2}+\frac{1}{2}
\gamma^{-1}\alpha^{-m/(r+m)}\right)^{-(r+m)}.
\end{equation}
Let $\delta=\delta(\alpha):=\left(1+\gamma^{-1}\alpha^{-m/(r+m)}\right)^{r+m}$.
The inequality (\ref{equ:convex1}) can be written
\[\forall t\notin K_\alpha(g)\quad  g(t)\leq 2^{-r-m}\delta g(t/2).
\]
Therefore
\[\E F(Y)1_{K_\alpha(g)^c}(Y)\leq
2^{-r-m}\int_{K_\alpha(g)^c} F(t)g\left(\frac{t}{2}\right)\delta\, dt
\leq 2^{-r}\int_{\R^m} F(2t)g(t)\delta\, dt
\]
and from the assumption on $F$, we get
\[\E F(Y)1_{K_\alpha(g)^c}(Y)
\leq 2^{p-r}\delta
\E F(Y).\]
We conclude that if  $2^{p-r}\delta<1$ then
\[ \E F(Y)\leq (1-2^{p-r}\delta)^{-1}\E \left(F(Y)1_{K_\alpha(g)}(Y)\right). \]
Let
\[\alpha _0=\left( \left(   2^{\frac{r-p}{2(r+m)}}-1 \right) \gamma\right)^{-\frac{r+m}{m}},\]
so that
$\delta _0:= \delta(\alpha_0)=2^{\frac{r-p}{2}}$, then
$(1-2^{p-r}\delta_0)^{-1}= (1-2^{\frac{p-r}{2}})^{-1}\leq 1+\frac{c}{r-p}$ and
$$\alpha_0\leq\alpha=\left(c\frac{(r+m)^2}{(r-p)(r-1)}\right)^{\frac{r+m}{m}},$$
where $c>0$ is a universal constant. This concludes the proof of (\ref{equ:restriction}).

Clearly $\gamma^{-1} \alpha^{-m/(r+m)}\leq \gamma^{-1} \alpha_0^{-m/(r+m)}<1$
and recall that $\gamma f(0)\leq \min f$. We deduce that
$f(0)\leq\alpha^{\frac{m}{r+m}} \min f$ and thus $0\in K_\alpha(g)$.
\qed

\medskip
\noindent{\bf Remark. } An interesting setting for the previous lemma is when $r$ is away
from~$1$, for instance $r\geq 2$,  $r$ and $m$ are comparable, and $p$
is proportional to $r$. In this case $\gamma$ is bounded by a constant,
$c(p,r)$ explodes only when $p\to r$,
and $\alpha$ depends only on the ratio $r/p$.

\bigskip
\noindent{\bf Proof  of Theorem \ref{theorem: assumption  convex}. }
Let $1\leq p<r$ and $ m=\lceil p \rceil$. Let $A:E\to \R^m$ be a linear mapping and
$Y=AX$ be a centered non-degenerated $(-1/r)$-concave random vector.
By Borell's result \cite{Bo1,Bo2}, there exists a positive convex
function $f$ such that the distribution of $Y$ has a density of the
form $g=f^{-(r+m)}$.

We apply Lemma~\ref{lemma: restriction} and use the notation of that lemma.
Because the class of $(-1/r)$-concave measures
 increases as the parameter $r$ decreases, we may assume that $r\leq 2p$
(note that $\lambda(p,2p)\sim \lambda(p,r)$ for $r>2p$, so we do not loose control
of the constant assuming that $r\leq 2p$). Thus
 $1\leq p\leq m$ and $r\leq 2m$. We deduce that the parameter  $\alpha$ from
Lemma~\ref{lemma: restriction}  satisfies
 $$\alpha\leq c\left(\frac{r}{r-1}\cdot\frac{r}{r-p}\right)^3, $$
 where $c$ is a numerical constant.

Now note that because $g^{-1/(r+m)}$ is convex,
 $K = K_\alpha(g)$ is a convex body and from Lemma~\ref{lemma: restriction},
it contains $0$. Let $\|\cdot\|$ be its Minkowski functional.

 We have
\[
   1 \geq \PP(Y \in K)=\int_{K}g \geq
   \alpha^{-m}\|g\|_\infty\text{vol}(K),
\]
so that
\[
  \PP(\| Y\|\leq 1/(2\alpha))=\int_{K/2\alpha}g \leq
  \|g\|_\infty(2\alpha)^{-m}\text{vol}(K) \leq 2^{-m} \leq 1/2,
\]
and therefore
$$
  \E \|Y\|\geq \frac{1}{2\alpha}\PP (\|Y\|>  1/(2\alpha))\geq  1/(4\alpha).
$$

Let $F(t)=\|t\|^p$ for $t\in\R^m$. Thus $F(2t)=2^pF(t)$ and, since $p<r$,
$\E F(Y)$ is finite.
Hence $F$ satisfies the assumption of Lemma~\ref{lemma: restriction}.
Therefore for $c(p,r)=1+c/(r-p)$
\[
  \E\| Y\|^p\leq c(p,r)\  \E \left( \|Y\|^p 1_{K}(Y)\right) \leq c(p,r).
\]
We conclude that
$$
  (\E\|Y\|^{p})^{1/p} /\E \|Y\|\leq 4\alpha\,  c(p,r)^{1/p}\leq
  c\left(\frac{r}{r-1}\right)^3\left(\frac{r}{r-p}\right)^4
$$
 for some numerical constant $c$.
\qed

\medskip

Another application of Lemma~\ref{lemma: restriction}, which will be used later,
 is the following lemma.

\begin{lemma}
\label{lemma: concave2}
Let  $1\leq p <r$ and $m=\lceil p\rceil$.
Let $Y\in\R^m$ be a centered $(-1/r)$-concave random
vector
with density $g$. There exists a universal constant $c$, such that
$0\in K_\alpha(g)$ and
\begin{equation}
\label{equ:restriction2}
\left(\E |\langle Y,t\rangle|^p\right)^{1/p}=
 \left(\int_{\R^m}  |\langle x,t\rangle|^p g(x)dx\right)^{1/p}
    \leq C_3(p, r)\,
    \max_{x\in K_\alpha(g)}  |\langle x,t\rangle|   ,
\end{equation}
where  $\alpha=c \, \left( \frac{r^2}{(r-p)(r-1)}\right)^3$,
$C_3(p, r)= \left(1+ \frac{c}{r-p}\right)^{1/p}$, and $c>0$ is a
universal constant.
\end{lemma}

\noindent
{\bf Proof. } Repeating the above proof with the function
$F(t)=|\langle x,t\rangle|^p$ we obtain that $0\in K_\alpha(g)$ and
$$
   \left(\int_{\R^m}  |\langle x,t\rangle|^p  g(x)dx\right)^{1/p}
 \leq \left(1+ \frac{c}{r-p}\right)^{1/p}\,\left(\int_{K_\alpha(g)}
    |\langle x,t\rangle|^p g(x)dx
    \right)^{1/p} .
$$
Clearly
$$
  \left(\int_{K_\alpha(g)}  |\langle x,t\rangle|^p  g(x)dx
  \right)^{1/p} \leq
   \max_{x\in K_\alpha(g)}  |\langle x,t\rangle| \left(\int_{\R^m}
     g(x)dx
    \right)^{1/p} = \max_{x\in K_\alpha(g)}  |\langle x,t\rangle| ,
$$
which implies the result.
\qed

\section{Small ball probability estimates}

The following result was proved in~\cite{Pa2}.
\begin{theorem}
\label{thm:pao3}
Let $X$ be a centered  log-concave random vector in a finite dimensional Euclidean space.
For every $\varepsilon\in (0,c')$ one has
$$
   \PP\left(|X|\leq \varepsilon (\E |X|^2)^{1/2}\right)\leq
   \varepsilon^{c(\E |X|^2)^{1/2}/{\sigma_2(X)}},
$$
where $c,c'>0$ are universal positive constants.
\end{theorem}

In this section we generalize this result to the setting of convex
distributions.
We first establish a lower bound for the negative moment of the
Euclidean norm of a convex random vector.

\begin{theorem}
\label{theorem:negative moment}
Let $r>1$ and let $X$ be a centered $n$-dimensional  $(-1/r)$-concave random vector. Assume $1\leq p < \min\{r, n/2\}$.
Then
$$
   \left( \E |X|^{-p}\right)^{-1/p} \geq  C_4 (p, r) \,  (\E |X| - C \sigma_p(X)),
$$
where
$$
 C_4(p, r) =c \, \left( \frac{r^2}{(r-p)(r-1)}\right)^{-3}
 \left(1+ \frac{c}{r-p}\right)^{-1/p}
$$ and  $c$, $C$ are absolute positive constants.
Moreover, if $0<p<1$ then
$$
  \left( \E |X|^{-p}\right)^{-1/p} \geq c_0 \,(1-p)  \, \frac{r-1}{r^2} \,  \E |X| ,
$$
where $c_0$ is an absolute positive constant.
\end{theorem}

\medskip

{}From Markov's inequality we deduce a small ball  probability
estimates for convex measures.

\begin{theorem}
\label{theorem:small ball}
Let $n\geq 1$ and $r>1$. Let $X$ be a centered $n$-dimensional
$(-1/r)$-concave random vector. Assume $1\leq p < \min\{r, n/2\}$.
Then, for every $\varepsilon\in (0,1)$,
\[
   \PP\left( |X| \leq \varepsilon \E |X|\right)\leq
    \left(2 C_4^{-1}(p, r) \varepsilon\right)^{p},
\]
whenever $\E|X|\geq 2 C\sigma_p(X)$, where $c$, $C$ and $C_4(p, r)$ are the
constants from Theorem~\ref{theorem:negative moment}.
\end{theorem}

\medskip

\noindent{\bf Remark. }   Theorem~\ref{theorem:small ball}
implies Theorem~\ref{thm:pao3} proved in \cite{Pa2}. Indeed,
let  $p\geq 1$,  $r\geq \max\{3, 2p\}$ and $A:=(\E|X|^2)^{1/2}/ \sigma_2(X)$
(note that $A\leq \sqrt{n}$).
By Lemma~\ref{lemma:khinchin},
$$
 \sigma_p(X)\leq C_1(p, r)\sigma_1(X)\leq c_0 p\sigma_2(X) \quad
  \mbox{ and } \quad (\E|X|^2)^{1/2} \leq c_1 \E|X|.
$$
Thus
$
 \E|X|/\sigma_p(X) \geq c_2 A/p .
$
If $c_2A/(2C)\geq 1$, we chose $p=c_2A/(2C)$ and apply
Theorem~\ref{theorem:small ball}. Since
$$
 \E|X|/\sigma_p(X) \geq c_2 A/p \geq 2 C,
$$
Theorem~\ref{thm:pao3} follows. Now assume that
$A\leq 2C/c_2$.
Then Theorem~\ref{thm:pao3}
follows from Lemma~\ref{mediantwo} and
Lemma~\ref{medianone} (with $q=2$), which for a log-concave random
vector $X$ states that
$
  \left(\E |X|^2 \right) ^{1/2} \leq c\, \mathrm{Med}(|X|),
$
where $c$ is a numerical constant and $\mathrm{Med}(|X|)$ is a median of $|X|$.

\bigskip

We need the following result from \cite{KV} (Theorem~1.3 there).

\begin{theorem}
\label{theorem:LO}
Let $n\geq 1$ be an integer, $\|\cdot\|$ be a norm in $\R^n$ and $K$ be its unit ball.
Assume that $0<p\leq c_0 \, \left( \E\|G\|/\sigma \right)^2$
and $m=\lceil p \rceil$.
Then
$$
 \frac{ c\, \E \|G\| }{ \sqrt{n} } \leq \left( \int_{G_{n,m}}  \left(
  \mbox{\rm diam}(K\cap F) \right)^{m} d\mu(F)
     \right)^{-1/m}\leq \frac{ \E \|G\|}{c  \sqrt{n}} ,
$$
where $\mu = \mu _{n,m}$ and $c$ is an absolute positive constant.
\end{theorem}

\medskip

The proof of Theorem~\ref{theorem:negative moment}
is based on the following two lemmas.

\begin{lemma}
\label{lemone}
Let $m\leq n$, $\alpha >0$ and $X$ be a random vector in $\R^n$ with
density $g$.
Then,
$$
  (\E |X|^{-m})^{-1/m}\geq \frac{1}{\sqrt{2\pi} \alpha}
  (\E |G|^{-m})^{-1/m} \left(\int_{G_{n,m}} \left(\mbox{\rm
        vol}(K_{\alpha}(g_F))\right)^{-1}d\mu(F)
  \right)^{-1/m}.
$$
\end{lemma}

\medskip

\noindent
{\bf Proof. }
Integrating in polar coordinates (see \cite{Pa2}, Proposition 4.6), we
obtain
the following key formula
$$
  (\E |X|^{-m})^{-1/m}= (2\pi)^{-1/2}(\E |G|^{-m})^{-1/m}
      \left(\int_{G_{n,m}}
          g_F(0)d\mu(F) \right)^{-1/m}.
$$
Note that
$$
 1= \int_F g_F(x) dx \geq  \int_{K_{\alpha}(g_F)} g_F(x) dx \geq
 \alpha^{-m} \|g_F\|_{\infty } \mbox{vol}(K_{\alpha}(g_F)).
$$
This implies the result, since $g_F(0) \leq  \|g_F\|_{\infty }$.
\qed

\medskip

Below we will use the following notation.
For a random vector $X$ in $\R^n$, $p>0$, and $t\in \R^n$
we denote
$$
  \|t\|_p = \left(\E  |\langle X, t\rangle|^p \right)^{1/p}
$$
(note that it is the dual gauge of the so-called centroid body, which is
rather an $L_p$-norm than the $\ell _p$-norm).

\begin{lemma}
\label{lemtwo}
Let $1\leq p < r$ and $m=\lceil p \rceil$. Let $X$ be a centered
$(-1/r)$-concave random vector in $\R^n$ with density $g$. Let $K$
denote the unit ball of $\|\cdot\|_p$. Then for every $m$-dimensional
subspace $F\subset \R^n$ one has
$$
   (\text{\rm vol}(P_F K^{\circ}))^{1/m}\leq 4 C_3 (p, r)
   (\text{\rm vol}(K_{\alpha}(g_F)))^{1/m},
$$
where  $\alpha=c \, \left( \frac{r^2}{(r-p)(r-1)}\right)^3$,
$C_3(p, r)= \left(1+ \frac{c}{r-p}\right)^{1/p}$, and $c>0$
is a universal constant.
\end{lemma}

\medskip

\noindent
{\bf Proof. }
Applying  Lemma~\ref{lemma: concave2} to $Y=P_FX$, we obtain that for
every $t\in F$
$$
   ||t||_p\leq C_3 (p, r)  \,\max_{x\in K_\alpha(g_F)}  |\langle x,t\rangle|
$$
with $\alpha$ and $C_3 (p, r) $ given in  Lemma~\ref{lemma: concave2}.
Since for $t\in F$, $||t||_p=\max \langle x, t\rangle$, where the supremum is taken over
$x\in (K\cap F)^{\circ} = P_F K^{\circ}$, this is equivalent to
$$
   P_F K^\circ \subset C_3 (p, r)  \text{\rm conv}(K_{\alpha}(g_F)\cup -K_{\alpha}(g_F)).
$$
 Lemma~\ref{lemma: concave2} also claims that $0\in K_{\alpha}(g_F)$, thus
$$
   \text{\rm conv}(K_{\alpha}(g_F)\cup -K_{\alpha}(g_F))\subset
   K_{\alpha}(g_F)     -K_{\alpha}(g_F).
$$
By Rogers-Sheppard  inequality \cite{RS} we observe
$$
  (\text{\rm vol}( P_F K^\circ))^{1/m}\leq {2m\choose m}^{1/m}
  C_3 (p, r)   (\text{\rm vol}(K_{\alpha}(g_F)))^{1/m}.
$$
This implies the result.
\qed

\medskip

\noindent
{\bf Proof of Theorem~\ref{theorem:negative moment}. }
Recall that  $c_1$, $c_2$, ... denote absolute positive
constants. Recall also that for a random vector $X$ in $\R^n$, $p>0$, and $t\in \R^n$
$$
  \|t\|_p = \left(\E  |\langle X, t\rangle|^p \right)^{1/p} .
$$
Given a norm $\|\, \cdot\, \|$ on $\R^n$
we define $\sigma$ by $\sigma=\sigma(\|\, \cdot\, \|) =\max_{|t|=1} \|t\|$.
In particular,
$$
  \sigma(\|\, \cdot\, \|_p) = \sigma _p (X).
$$
Finally, let   $K$ denote the unit ball of  $\|\, \cdot\, \|_p$.

We assume that $X$ is non-degenerate in $\R^n$ and
let $m=\lceil p \rceil$.
Without loss of generality we assume that
$$
  \E |X|\geq C   \sigma_p(X) ,
$$
where $C$ is a large enough  absolute constant.

As in (\ref{equ:gaussian moment}), since $p\leq m \leq 2p$, we have
$$
 \E |X|\leq (\E |X|^p)^{1/p}= \gamma_p^{-1} (\E\|G\|_{p}^p)^{1/p}
  \leq  \gamma_p^{-1}( \E\|G\|_{p}+ c_1  \gamma_p \sigma_p(X) )
$$
$$
 \leq c_2 (\E\|G\|_{p}/\sqrt{m} + \sigma_p(X) ) .
$$
Hence
\begin{equation}\label{dvorlevel}
  \E\|G\|_{p} \geq \sqrt{p}  c_2^{-1}  \left(\E |X| - c_2  \sigma_p(X) \right)
 \geq  \sqrt{p} (c_2)^{-1} (C-c_2)   \sigma_p(X).
\end{equation}
This implies that for sufficiently large $C$ we have
$m\leq 2p\leq c_0 (\E \|G\|_p/\sigma _p(X))^2$, where $c_0$ is the constant
from Theorem~\ref{theorem:LO}.

Note that
$ ( \E |G|^{-p})^{-1/p} \geq ( \E |G|^{-m})^{-1/m} \geq c_3 \sqrt{n}$
(the second inequality is well known for $m\leq n/2$ and can be directly computed).
Combining  Lemmas \ref{lemone} and \ref{lemtwo}, we obtain
$$
  (\E |X|^{-m})^{-1/m}\geq \frac{c_4 \sqrt{n}}{\alpha\, C_3 (p, r)} \,
  \left(\int_{G_{n,m}} (\text{\rm vol}(P_F K^{\circ}))^{-1}d\mu(F)
  \right)^{-1/m},
$$
with $\alpha$ and $C_3(p, r)$ as in Lemma~\ref{lemtwo}.

Now note that
$P_F K^{\circ} = (K\cap F)^{\circ} \supset (\mbox{diam}(K\cap F))^{-1} B_2^n \cap F$.
Therefore $1/\mbox{vol}(P_F K^{\circ}) \leq (c_5\sqrt{m}\, \mbox{diam}(K\cap F))^{m}$.
Applying  Theorem~\ref{theorem:LO}, we obtain
$$
  (\E |X|^{-m})^{-1/m}\geq \frac{c_6 }{\alpha\, C_3 (p, r) \sqrt{m}} \,
    \E \|G\|_p .
$$
Applying the first inequality from (\ref{dvorlevel}), we obtain the desired result.

The ``Moreover" part is an immediate corollary of Lemmas~\ref{medianone} (with $q=1$)
 and \ref{mediantwo}.
\qed

\bigskip

\noindent{\bf Conjecture.}  We conjecture that for
convex distributions a similar  thin shell property holds as for log-concave
distribution: if $X$ is an isotropic  $(-1/r)$-concave random vector in $\R^n$ with $r>2$, then
$$
  \forall t\in (0,1)\quad \PP\left(\big| |X| - \E|X| \big|\geq t\sqrt n\right)
  \to  0 .
$$
as $n$ tends to $\infty$.
See \cite{GM} for recent work in the log-concave setting.

\section{Appendix}

There is a vast literature on  inequalities of integrals related to concave functions. Some of the
following lemmas may be known but we did not find any reference. Their proof use classical methods for
demonstrating integral inequalities involving concave functions (see \cite{Bo3} and \cite{MP}). The
first lemma is a mirror image for negative moments of a result from \cite{MOP} valid for positive moments.

\begin{lemma}
\label{lemma:Gincreasing}
Let $s,m,\beta\in\R$ such that $\beta>m+1>0$ and $s>0$. Let $\varphi$ be a non-negative concave function on
$[s,+\infty)$. Then
\[G(\beta)= \frac{\int_s^\infty \varphi^m (x) x^{-\beta}\,dx}{s^{m-\beta+1}B(m+1,\beta-m-1)}\]
is an increasing function of $\beta$ on $(m+1,\infty)$.
Here $B(u,v)=\int_0^1(1-t)^{u-1}t^{v-1}\, dt$ denotes the Beta function.
\end{lemma}

\noindent{\bf Proof. }
Let $\beta > m+1$.
Consider the function
\[ H(t)=\int_s^t \varphi^m(x) x^{-\beta}\,dx-
\int_s^t a^m(x-s)^m x^{-\beta}\,dx\]
for $t\geq s$, where $a$ is chosen so that $H(\infty)=0$.
Note that $H'$, the derivative of $H$, has the same sign as $(\varphi(x)/(x-s))^m-a^m$.
Since $\varphi(x)/(x-s)$ is decreasing on $(s,+\infty)$, we deduce
that $H$ is first increasing and then decreasing. Since $H(s)=H(\infty)=0$
we conclude that $H$ is non-negative. This means that for every $t\geq s$,
\begin{equation}
\label{equ:H}
\int_s^t \varphi^m(x) x^{-\beta}\,dx\geq \int_s^t a^m(x-s)^m x^{-\beta}\,dx.
\end{equation}
Now, note that for any $\beta'> \beta$ and any non-negative function $F$, we have
by Fubini's theorem,
\[\int_s^\infty F(x) x^{-\beta'}\,dx=\int_s^\infty (\beta'-\beta) t^{-\beta'+\beta-1}\left(\int_s^t F(x) x^{-\beta}\,dx\right)dt.\]
Using (\ref{equ:H}) and applying this relation to $F=\varphi^m$ and then to $F(x)=a^m(x-s)^m$,
we get that
\[
  \int_s^\infty \varphi^m(x) x^{-\beta'}\,dx\geq a^m\int_s^\infty (x-s)^m x^{-\beta'}\,dx.
\]
{}From the definition of $a$, we conclude that
\[{\int_s^\infty \varphi^m(x) x^{-\beta}\,dx}\Big/{\int_s^\infty (x-s)^m x^{-\beta}\,dx}\]
is an increasing function of $\beta$ on $(m+1,\infty)$. The conclusion follows from the computation
of $\int_s^\infty (x-s)^m
x^{-\beta}\,dx=s^{m-\beta+1}B(m+1,\beta-m-1)$.\qed

\begin{lemma}
\label{corollary:g(0)}
Let $m\geq 1$ be an integer. Let $g$ be the density  of a probability on $\R^m$ of the form
$g=f^{-\beta}$ with $f$ positive convex on $\R^m$ and $\beta> m+1$.
If $\int xg(x)\,dx=0$, then
$$
g(0)\geq \left(\frac{\beta-m-1}{\beta-1}\right)^\beta \|g\|_\infty.
$$
\end{lemma}

\noindent{\bf Proof. } Since $f$ is convex it follows from Jensen's inequality that
\[ f(0)=f\left(\int xg(x)\,dx\right)\leq \int f^{-\beta +1}(x)\, dx
=\int_s^\infty (\beta -1)h(s) s^{-\beta}\, ds,\]
where $s=\min f=\|g\|_\infty^{-1/\beta}$ and $h(t)=\text{vol}\{f\leq t\}$ denotes the Lebesgue measure
of $\{f\leq t\}$. From the convexity of $f$ and from the Brunn-Minkowski inequality, $\varphi=h^{1/m}$ is
concave.  Thus, using the notation of Lemma~\ref{lemma:Gincreasing},
\[f(0)\leq (\beta-1)s^{m-\beta+1}B(m+1,\beta-m-1)G(\beta).
\]
Now observe that $\int f^{-\beta}=\int_s^\infty \beta \varphi^m(x) x^{-\beta-1}\,dx=1$ and
therefore, by Lemma~\ref{lemma:Gincreasing},
\[G(\beta)\leq G(\beta+1)=\left(\beta s^{m-\beta}B(m+1,\beta-m)\right)^{-1}.\]
The conclusion follows from combining the last two inequalities.\qed

\medskip
\noindent {\bf Remark. } When $\beta\to\infty$, which corresponds to a log-concave density,
 we recover the inequality from \cite{Fra} saying that
 $g(0)\geq e^{-m} \|g\|_\infty$.

\medskip

The next lemma is a well known result  of Borell (\cite{Bo1}) stated
in a way that fits our needs
and  stresses the dependence on the parameter of concavity.

\begin{lemma}
\label{lemma:khinchin} Let $r>1$ and $X$ be a $(-1/r)$-concave random vector in $\R^m$.
Then for any semi-norm $\|\,\cdot \, \|$ and any $t\geq 1$, one has
\[
   \PP(\|X\|\ \geq 3t \E  \|X\|)  \leq  \left(1+ \frac{ t }{ 3 r} \right)^{-r}.
\]
As a consequence, for every $1\leq p<r$,
\[
    (\E \|X\|^p)^{1/p}\leq C_1(p, r) \,\E \|X\|,
\]
where
$C_1(p,r)=cp$ for $r>p+1$, $C_1(p,r)=\frac{c\, r}{(r-p)^{1/p}}$
otherwise and
$c$  is a universal constant.
\end{lemma}

\noindent{\bf Proof. } Denote $\theta:=\PP\left(\|X\| \leq 3\E  \|X\|\right)$.
Assume that $\theta<1$ (otherwise we are done).
{}From Markov's inequality,
\[\theta=
 1-\PP(\|X\|\, > 3 \E  \|X\|) \geq 2/3.\]
The subset $B=\{x\in\R^m: \|x\| \leq 3\E  \|X\|\}$ is symmetric and convex.
 From Lemma 3.1 in \cite{Bo1}, for every $t\geq 1$, one has
\[\PP\left(\|X\|\ \geq 3 t\E  \|X\|\right)\leq
\left(     \frac{t+1}{2}\left( (1-\theta)^{-1/r}-\theta^{-1/r}\right) +\theta^{-1/r}       \right)^{-r}.\]
Thus,
\[\PP\left(\|X\|\ \geq 3 t\E  \|X\|\right)\leq
\theta\left(   1+ \frac{1}{2r}\log\frac{\theta}{1-\theta} +\frac{t}{2r}
\log\frac{\theta}{1-\theta}  \right)^{-r} .\]
We deduce that for every $t\geq 1$,
\[\PP\left(\|X\|\ \geq 3 t\E  \|X\|\right)\leq
\theta\left(   1+ \frac{t}{2r}\log\frac{\theta}{1-\theta}  \right)^{-r}
 \leq\left(   1+ \frac{t}{3r}\right)^{-r}.\]
Integrating, we get
\begin{align*}
  \E \|X\|^p\big/(3\E \|X\|)^p
 &=
   \int_0^\infty pt^{p-1}\PP\left(\|X\|\ \geq 3 t\E  \|X\|\right) dt     \\
 &\leq 1+\int_1^\infty pt^{p-1}\left(   1+ \frac{t}{3r}\right)^{-r} \,dt \\
 &\leq 1+(3r)^p pB(p, r-p)                        \\
  &= 1+(3r)^p\Gamma(p+1)\Gamma(r-p)/\Gamma(r).
\end{align*}
Now, if $r>p+1$ then, by Stirling's formula,
$$
  \left( \Gamma(p+1)\Gamma(r-p)/\Gamma(r) \right)^{1/p} \sim \frac{p}{r},
$$
and if $r\leq p+1$ then
$$
  \left( \Gamma(p+1)\Gamma(r-p)/\Gamma(r) \right)^{1/p} \sim
  \left( \Gamma(r-p) \right)^{1/p} \sim \left( \frac{1}{r-p} \right)^{1/p} .
$$
This completes the proof.
\qed

\medskip

The following, a stronger variant of Borell's lemma, allows to compare
the expectation of a random variable $\|X\|$  and a median $\mathrm{Med}(\|X\|)$.
It was proved in \cite{Bobkov} (Theorem~1.1, see also the discussion
following Theorem 5.2 for the behavior of the corresponding constant).
It was also implicitly
proved in \cite{G} (see inequality (4) in \cite{Fr}). The second part of the
lemma follows by integration.

\begin{lemma}
\label{medianone} Let $r>1$ and $X$ be a $(-1/r)$-concave random vector in $\R^m$.
Then for any semi-norm $\|\,\cdot \, \|$ and any $t\geq 1$, one has
\[
   \PP\left(\|X\|\ \geq t \mathrm{Med}( \|X\|)  \right)  \leq   (C_0 r)^r t^{-r},
\]
where $C_0$ is an absolute positive constant.
As a consequence, for every $r>q\geq 1$  one has
$$
  \left(\E \|X\|^q \right) ^{1/q} \leq C r \left(\frac{r}{r-q}\right)^{1/q}\,
  \mathrm{Med}(\|X\|),
$$
where $C$  is an absolute positive constant.
\end{lemma}

The following lemma is Corollary~9 from \cite{Fr} (as before, the second
part follows by  integration).

\begin{lemma}
\label{mediantwo} Let $r>1$ and $X$ be a $(-1/r)$-concave random vector in $\R^m$.
Then for any semi-norm $\|\,\cdot\, \|$ and any $\eps\in (0,1)$, one has
\[
   \PP\left(\|X\|\ \leq \eps \mathrm{Med}( \|X\|)  \right)  \leq  C_0 \eps,
\]
where $C_0$ is an absolute positive constant.
As a consequence, for every $p\in (0,1)$,
\[
    \left( \E \|X\|^{-p}\right)^{-1/p}  \geq c (1-p) \mathrm{Med}(\|X\|),
\]
where $c$  is an absolute positive constant.
\end{lemma}

\address

\end{document}